\documentclass[a4paper, 11pt]{article} 
\pdfoutput=1

\usepackage{graphicx} 
\usepackage{hyperref}
\usepackage{amsmath}
\usepackage{amssymb}
\usepackage{amsthm}
\newtheorem*{conjecture}{Conjecure}
\newcounter{savefootnote}
\usepackage[T1]{fontenc} 

\makeatletter
\renewcommand\@biblabel[1]{\textbf{#1.}} 
\renewcommand{\@listI}{\itemsep=0pt} 

\renewcommand{\maketitle}{ 
\begin{flushright} 
{\LARGE\@title} 

\vspace{50pt} 

{\large\@author Dimitris Vartziotis 
	\footnote{TWT GmbH Science \& Innovation, Department for Mathematical Research, Ernsthaldenstra{\ss}e 17, 70565 Stuttgart, Germany}\setcounter{savefootnote}{\value{footnote}}\footnote{NIKI Ltd. Digital Engineering, Research Center, 205 Ethnikis Antistasis Street, 45500 Katsika, Ioannina, Greece}
	\footnote{Corresponding author. E-mail address: dimitris.vartziotis@nikitec.gr},
	Juri Merger\setcounter{footnote}{\value{savefootnote}}\footnotemark[\value{footnote}]} 
\\\@date 

\vspace{40pt} 
\end{flushright}
}

\title{\textbf{Contributions to the study of the non-trivial roots of the Riemann zeta-function}\\ 
} 

\author{\textsc{} 
\\{\textit{}}} 

\date{\today} 

\begin{document}

\maketitle 


\renewcommand{\abstractname}{Summary} 

\begin{abstract}
This work contributes to the study of the non-trivial roots of the Riemann zeta function.
In view of the Hilbert-P\'{o}lya conjecture a series of self-adjoint operators on a Hilbert space is constructed whose eigenvalues approximate these roots.
\end{abstract}

\hspace*{3,6mm}\textit{Keywords:} Riemann zeta function, Fourier analysis, Hilbert-P\'{o}lya conjecture, Montgomery's pair correlation conjecture 

\vspace{30pt} 


\section{Introduction}
The nature of prime numbers fascinates mathematicians since many centuries. 
Its various connections to almost all fields for mathematics shows their special role.
E.g. prime numbers can induce fractal structures, see \cite{vartziotis2017fractal, vartziotis2018fractal}, or have relations to polygonal transformations \cite{vartziotis2009classification}.
Central to this paper is the Riemann hypothesis (RH) saying that all non-trivial roots of the Riemann zeta function lie on the critical line with real part $\frac{1}{2}$.

Many people also tried to underpin (RH) by numerical computations. Therefore, there are large data sets available on the internet containing the first roots to a certain accuracy \cite{lmfdb}.
See also the corresponding literature \cite{odlyzko200110,odlyzko1988fast}.

In this work we set the connection between Riemann roots and hermitian, cyclic matrices also used in the work of \cite{VARTZIOTIS2010945}.

Moreover, it seems that the values of the Riemann zeta-function can be modelled by random matrices, see \cite{conrey2005integral}.

A comprehensive study of the value distributions of $L$-functions (a generalization of the Riemann zeta function) can be found in \cite{steuding2007value}.

\newpage
\section{Fourier transform}
\begin{figure}[]
	\centering
	\includegraphics[width=\textwidth]{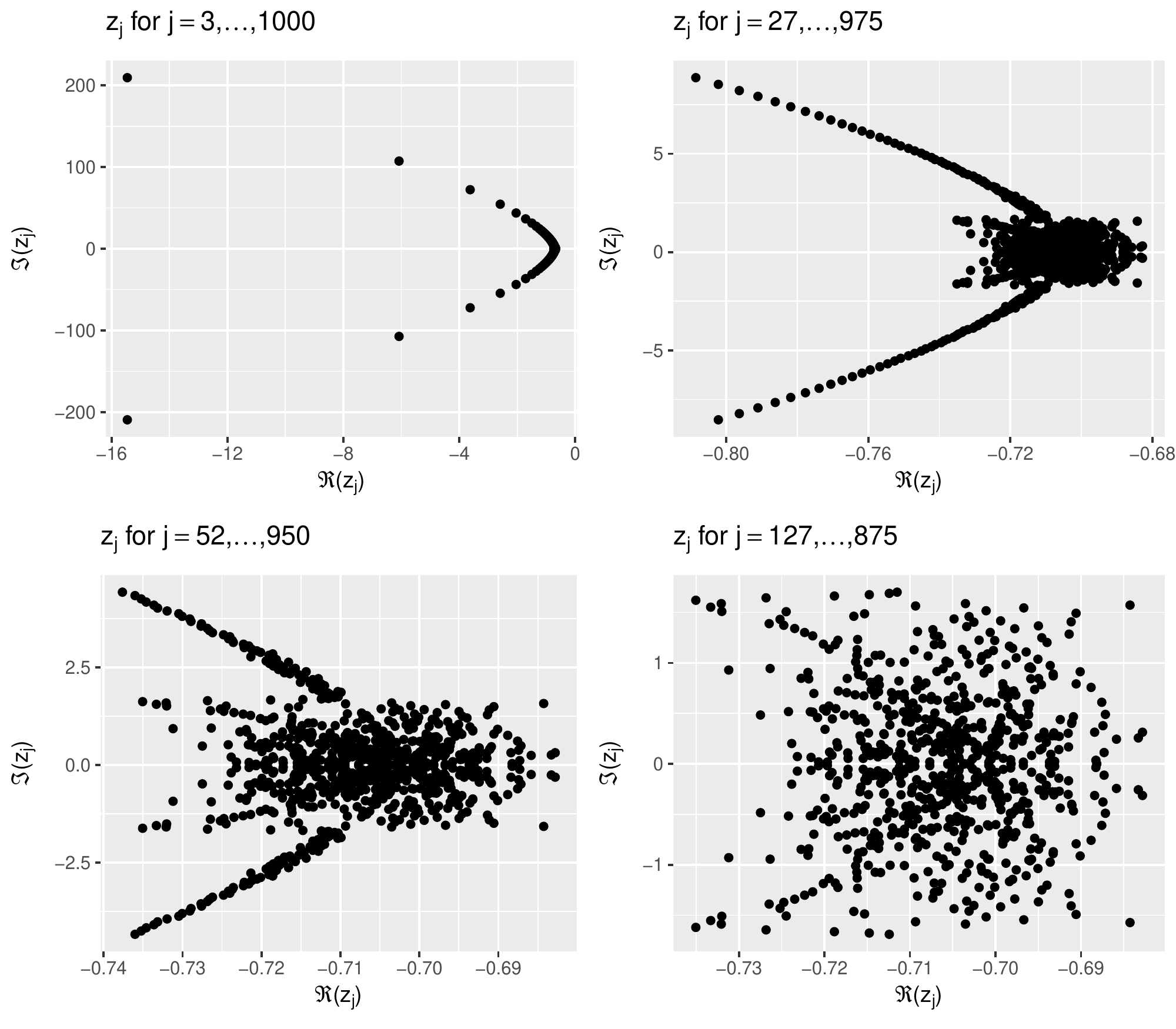}
	\caption{Different zoom levels for the distribution of Fourier coefficients for the first 1000 roots $\gamma_j$ $(j = 1,\ldots,1000)$.}
	\label{fig:fourier_coef}
\end{figure}

In this work we analyze the discrete Fourier transform of the non-trivial roots $\gamma_n$ $(n \in \mathbb{N})$ of the Riemann zeta function. For each $n\in\mathbb{N}$ they are defined as follows
\begin{equation}
z_k^{(n)} = \sum_{j = 1}^{n} \frac{1}{n} \gamma_j \exp\left(\frac{2 \pi i}{n} (j-1)(k-1)\right)
\label{eq:def_of_zn}
\end{equation}
Note, that the first Fourier coefficient $z_1^{(n)} = \sum_{j = 1}^{n} \frac{1}{n} \gamma_j$ is the mean of the first $n$ roots $\gamma_j$. 
The distribution of remaining Fourier coefficient can be seen in \autoref{fig:fourier_coef}.
It seems that the $z_j$ follow globally a curve that looks like a left-opened parabola.
Whereas most of the values are clustered at the vertex of it.

\begin{figure}[]
	\centering
	\includegraphics[width=\textwidth]{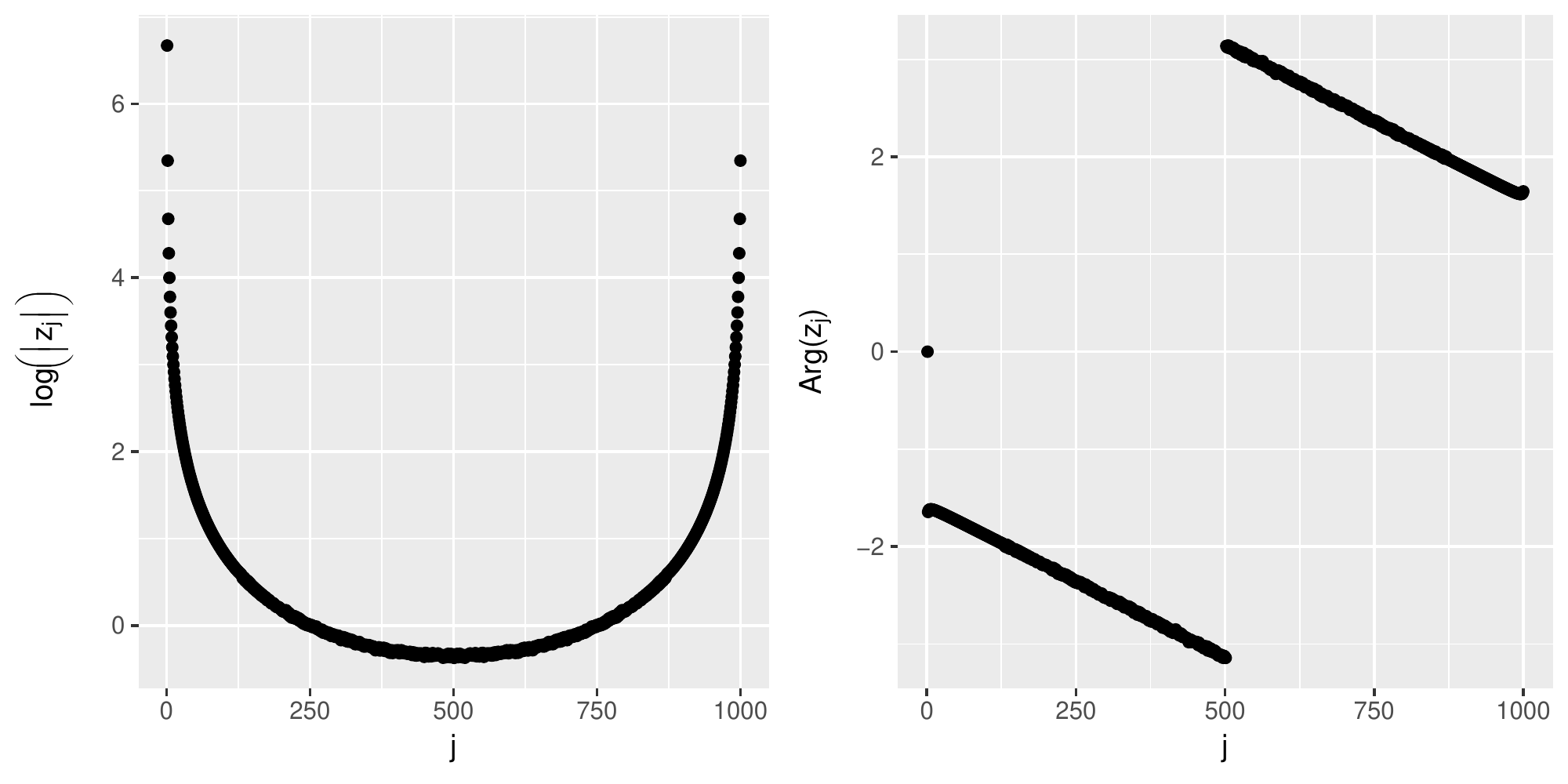}
	\caption{Modulus and argument of $z_j^{(1000)}$ for $j = 1,\ldots, 1000$.}
	\label{fig:fourier_coef2}
\end{figure}

Moreover, the modulus and argument of the Fourier coefficients are viewed in \autoref{fig:fourier_coef2}.
Surprisingly, the argument is almost linear in $j$.

The Fourier coefficients $z_k^{(n)}$ encode all the characteristics of the series of roots $\gamma_n$.
Therefore, we can reconstruct them by the following formula
\begin{equation}
\gamma_k = \sum_{j = 1}^{n} z_j^{(n)} \exp\left(-\frac{2 \pi i}{n} (j-1)(k-1)\right).
\label{eq:reconstruct_roots}
\end{equation}

\begin{figure}[]
	\centering
	\includegraphics[width=\textwidth]{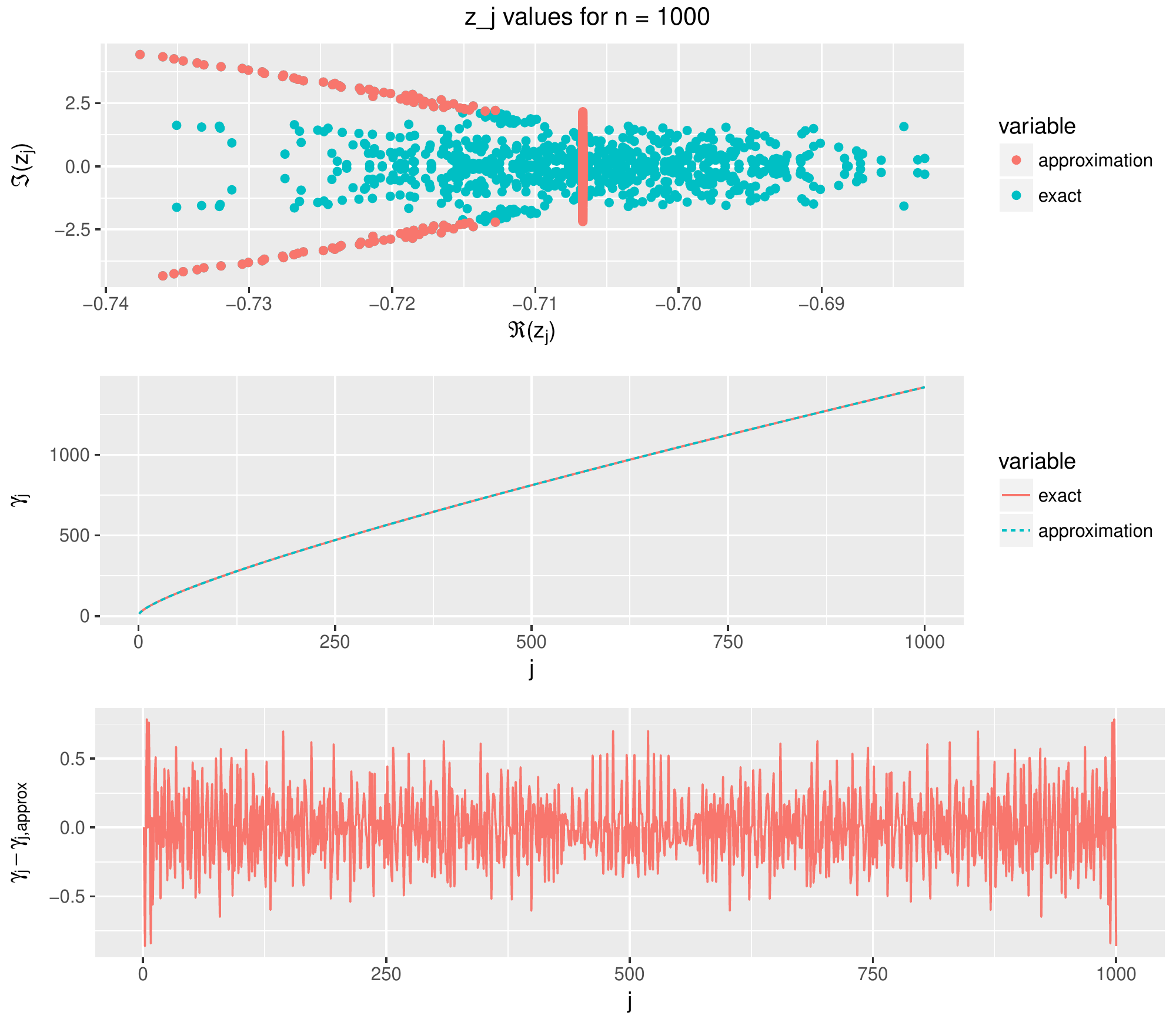}
	\caption{For the 80 \% intermediate Fourier coefficients ($n = 1000$) the real part was averaged (upper plot) and compared the approximation of the roots $\gamma_{j, approx}$ reconstructed from them to the exact values ${\gamma}_j$ (absolute values: middle plot; differences: lower plot).}
	\label{fig:constant_real_part}
\end{figure}

It is known that the global structure of the sequence $\gamma_1, \ldots, \gamma_n$ is encoded in the first and last few Fourier coefficients $z_j^{(n)}$ with $j \approx 1$ or $j \approx n$, whereas the local differences of subsequent roots $\gamma_j, \ldots, \gamma_{j+k}$ with $k << n$ are encoded in the intermediate values $z_j{(n)}$ for $j \approx \frac{n}{2}$.
Therefore, knowing the exact real part of most Fourier coefficients is not necessary to reconstruct the asymptotic behavior of the series $\gamma_n$.
See \autoref{fig:constant_real_part}, where we exchange the real part of 80 \% of the Fourier coefficients by its mean and can still reconstruct the global curve described by the roots $\gamma_n$.

\section{Montgomery pair correlation conjecture}
According to \cite{montgomery1973pair}, we can approximate under the Riemann hypothesis the imaginary part $\gamma_n$ of the $n$-th root $\frac{1}{2} + \gamma_n i$ of the zeta function by the following recursion formula:
\begin{equation}
\gamma_{n+1} = \gamma_n + \frac{2 \pi}{\log \frac{\gamma_n}{2 \pi}}
\label{eq:montgomery_approx}
\end{equation}
The approximation according to \eqref{eq:montgomery_approx} is very good and the relative error goes to zero.

In the following we want to study not the roots $\gamma_n$ itself, but the discrete Fourier transform of them. 
Therefore, define for each integer $n$ the following sequence:
\begin{equation}
z_k^{(n)} = \sum_{j = 1}^{n} \frac{1}{n} \gamma_j \exp\left(\frac{2 \pi i}{n} (j-1)(k-1)\right)
\label{eq:def_of_zn}
\end{equation}

Moreover, the first value $z_1^{(n)}$ is the mean of the first $n$ roots. Therefore, we have
$$ z_1^{(n)} = \sum_{j = 1}^{n} \frac{1}{n} \gamma_j = \frac{n-1}{n}\sum_{j = 1}^{n-1} \frac{1}{n-1} \gamma_j + \frac{\gamma_n}{n} = \frac{n-1}{n} z_1^{(n-1)} + \frac{\gamma_n}{n}$$
and the $n$-th root is given by
\begin{equation}
\gamma_n = n z_1^{(n)} - (n-1)z_1^{(n-1)}.
\end{equation} 

To derive a recursive formula for $z$-values, we compute
\begin{align*}
n z_k^{(n)} 
=& \sum_{j = 1}^{n-1} \gamma_j \exp\left(\frac{2 \pi i}{n} (j-1)(k-1)\right) + \gamma_n \exp\left(\frac{2 \pi i}{n} (n-1)(k-1)\right)\\
=& \sum_{j = 1}^{n-1} \left[\sum_{l = 1}^{n-1} z_l^{(n-1)} \exp\left(-\frac{2 \pi i}{n-1} (l-1)(j-1)\right) \right] \exp\left(\frac{2 \pi i}{n} (j-1)(k-1)\right) \\
&+ \gamma_n \exp\left(-\frac{2 \pi i}{n}(k-1)\right) \\
=& \sum_{l = 1}^{n-1} z_l^{(n-1)} \left[\sum_{j = 1}^{n-1} \exp\left[2 \pi i \left(\frac{k-1}{n} - \frac{l-1}{n-1}\right)(j-1)\right]  \right]\\
&+ \left(\gamma_{n-1} + \frac{2 \pi}{\log \frac{\gamma_{n-1}}{2 \pi}}\right) \exp\left(-\frac{2 \pi i}{n}(k-1)\right) \\
=& \sum_{l = 1}^{n-1} z_l^{(n-1)} \left(1 + \xi_{k,l} + \xi_{k,l}^2 + \ldots + \xi_{k,l}^{n-2} \right) \\
&+ \left((n-1) z_1^{(n-1)} - (n-2)z_1^{(n-2)} + \frac{2 \pi}{\log \frac{(n-1) z_1^{(n-1)} - (n-2)z_1^{(n-2)}}{2 \pi}}\right) \\
& \times \exp\left(-\frac{2 \pi i}{n}(k-1)\right),
\end{align*}
where we define the $\frac{n(n-1)}{n(k-l) + 1 - k} $-th root of unity by
$$\xi_{k,l} :=  \exp\left(2 \pi i \left(\frac{k-1}{n} - \frac{l-1}{n-1}\right)\right). $$

Hence, we arrive at the following recursive formula
\begin{align}
z_k^{(n)} =
& \frac{1}{n}\sum_{l = 1}^{n-1} z_l^{(n-1)} \left(1 + \xi_{k,l} + \xi_{k,l}^2 + \ldots + \xi_{k,l}^{n-2} \right) \\
&+ \frac{1}{n}\left((n-1) z_1^{(n-1)} - (n-2)z_1^{(n-2)} + \frac{2 \pi}{\log \frac{(n-1) z_1^{(n-1)} - (n-2)z_1^{(n-2)}}{2 \pi}}\right) \\
& \times \exp\left(-\frac{2 \pi i}{n}(k-1)\right)
\label{eq:fourier_approx}
\end{align}

\begin{figure}[]
	\centering
	\includegraphics[width=\textwidth]{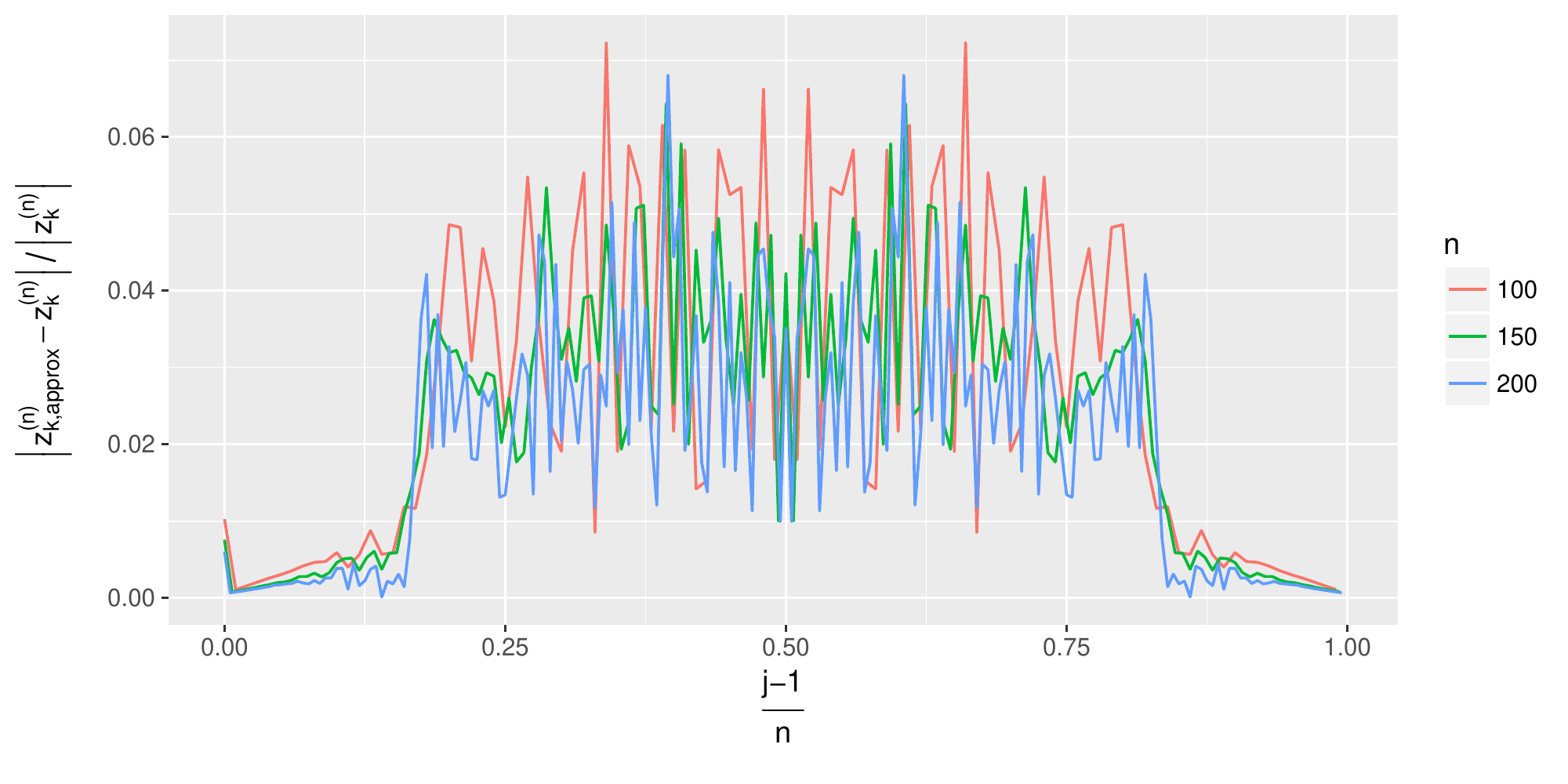}
	\caption{Relative difference of the exact $z$ values and the approximations according to \eqref{eq:fourier_approx}.}
	\label{fig:fourier_approx}
\end{figure}

Here, we observe also a decrease of the relative error as $n$ tends to infinity; see \autoref{fig:fourier_approx}.

By reconstructing the roots from the approximate Fourier coefficients we end up with an approximation for the non-trivial root of the Riemann zeta function. 
A comparison to the exact values shows the same results as \eqref{eq:montgomery_approx}.

\begin{figure}[]
	\centering
	\includegraphics[width=\textwidth]{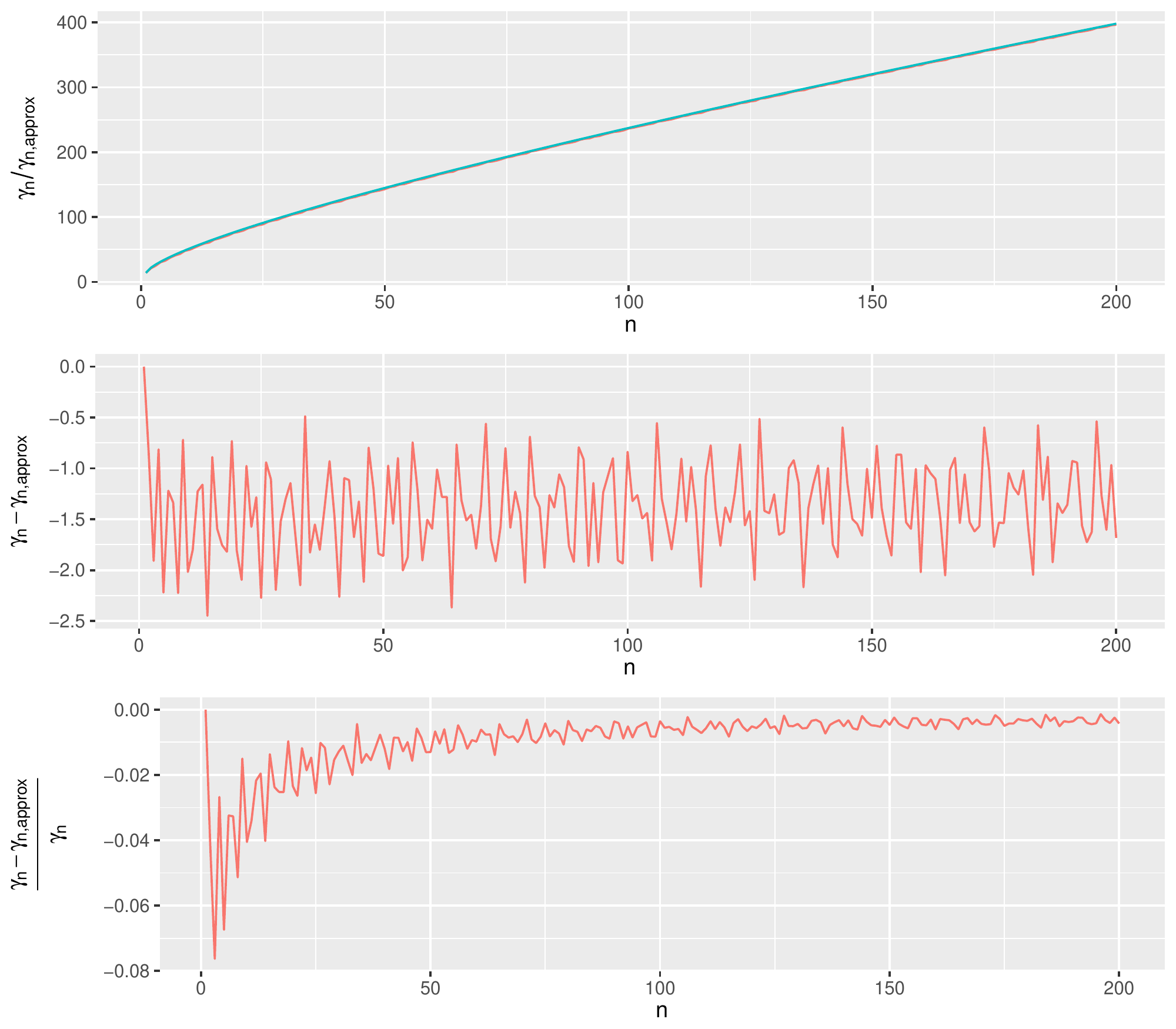}
	\caption{A comparison of the exact roots of the Riemann zeta function and the approximation by the recursion formula for the Fourier coefficients $z_j^{(n)}$.}
	\label{fig:matrix_approx}
\end{figure}

\begin{conjecture}
\textbf{Hilbert-Polya}:
If we write the non-trivial Riemann zeros as $\frac{1}{2} + i \gamma_n$, where $\Re(\frac{1}{2} + i \gamma_n) > 0$, then the numbers $\gamma_n$ correspond to the eigenvalues of an unbounded self-adjoint operator.
\end{conjecture}
	
Our aim is to approximate such an unbounded self-adjoint operator and show its convergence. Therefore, we define the following cyclic matrix
\begin{equation}
M_n = \left( z_{1+k-l}^{(n)} \right)_{k,l = 1,\ldots,n}
\label{eq:def_M}
\end{equation}
where $z_{k}^{(n)}$ is defined as above.
It is known that the normalized eigenvectors of a cyclic matrix are given the Fourier-polygons 
$$ v_j = \frac{1}{\sqrt{n}}\left(\exp\left(\frac{2 \pi i}{n}(k-1)(j-1)\right)\right)_{k = 1,\ldots,n} $$
and, therefore, its eigenvalues are
\begin{align*}
\lambda_j =& v_j^* M_n v_j = \sum_{k,l = 1}^n \frac{1}{n} z_{1+k-l}^{(n)} \exp\left(-\frac{2 \pi i}{n}(k-l)(j-1)\right) \\
=& \sum_{k = 1}^n z_{k}^{(n)} \exp\left(-\frac{2 \pi i}{n}(k-1)(j-1)\right) = \gamma_j
\end{align*} 
the imaginary part of the non-tivial roots of the zeta function.

\section{Construction of an self-adjoint operator}
Let $H$ be a complex, separable Hilbert space (e.g. $L^2(0,1)$) and $(e_j)_{j = 1,\ldots,\infty}$ an orthonormal basis of $H$. Then define the following sequence of operators
\begin{align*}
T_n : H &\rightarrow H \\
x &\mapsto \frac{1}{2} + i * P_n^{-1} (M_n P_n(x)),
\end{align*}
where $P_n:H \rightarrow \mathbb{R}^n$ is defined by $P_n(x) = \left(\left\langle x, e_j \right\rangle \right)_{j = 1,\ldots,n}$ and $P_n^{-1}:\mathbb{R}^n \rightarrow H$ is defined by $P_n^{-1}(x) = \sum_{j = 1}^{n} x_j e_j$. Then, $T_n$ is a self-adjoint operator whose eigenvalues are given by the first $n$ roots of the Riemann zeta function.
If the limit of $T_n$ as $n$ tends to infinity exists the Hilbert-Polya conjecture could be proven.

\section{Additional statistical analysis}
In the following we document some further analysis that we have done with respect to the Fourier coefficients of the Riemann root $\gamma_n$.

\subsection{Analysis of the Fourier coefficient cluster}
For all $n \in N$ the plot of the Fourier coefficients seems to be similar.
As can be seen in \autoref{fig:fourier_coef}, the global distribution of these values follow a smooth curve, whereas most of the intermediate values ($j \approx \frac{n}{2}$) form a cluster at the vertex of the left open curve.
Moreover, the real part of the coefficients seem to have a very high variance compared to the "smoother" regions for ($j \approx 1$ or $j \approx n$).
In the following, we try to separate those two different regions by finding an index $m \in \{1, \ldots, n\}$ such that the coefficients $z_m, \ldots, z_{n-m}$ form the cluster and the remaining values constitute the global smooth curve without the vertex.

\begin{figure}[]
	\centering
	\includegraphics[width=\textwidth]{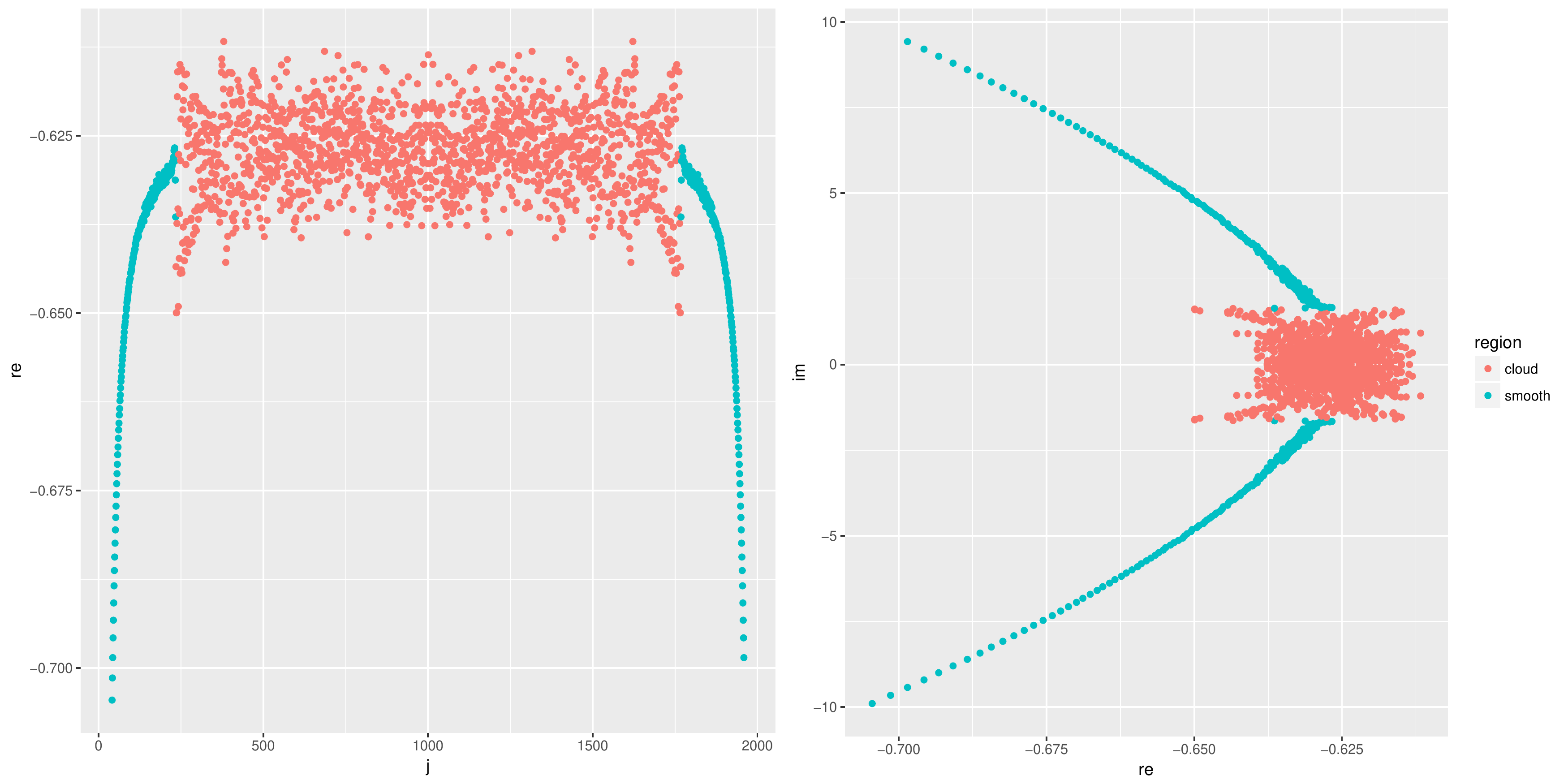}
	\caption{The results of the cloud index determination algorithm is shown from two perspectives. The left plot shows the real part of the Fourier coefficients this respect to its index $j$. The right plot shows the distributions of $z_j$ across the complex plane.}
	\label{fig:cloud_index}
\end{figure}

In order to achieve this goal we utilize the increasing variance of the real part.
In particular, the algorithm goes as follows.
First, we choose a window width $d << n$ and consider the real part of the coefficients for $j - d, \ldots, j + d$. 
Then, after normalizing to mean 0 and maximum 1 we apply a discrete Fourier transform in order to measure the smoothness of this section.
At last, we compute the weighted mean of the absolute values of these Fourier coefficients such that the first half has positive weights and the second half has negative.
This leads in total to a negative result provided the considered values in the window at hand have a high variance and are not smooth, such that the higher Fourier coefficients have a high absolute value.
In the end we take the first index where the weighted mean is negative as the beginning of the cloud. 

\begin{figure}[]
	\centering
	\includegraphics[width=\textwidth]{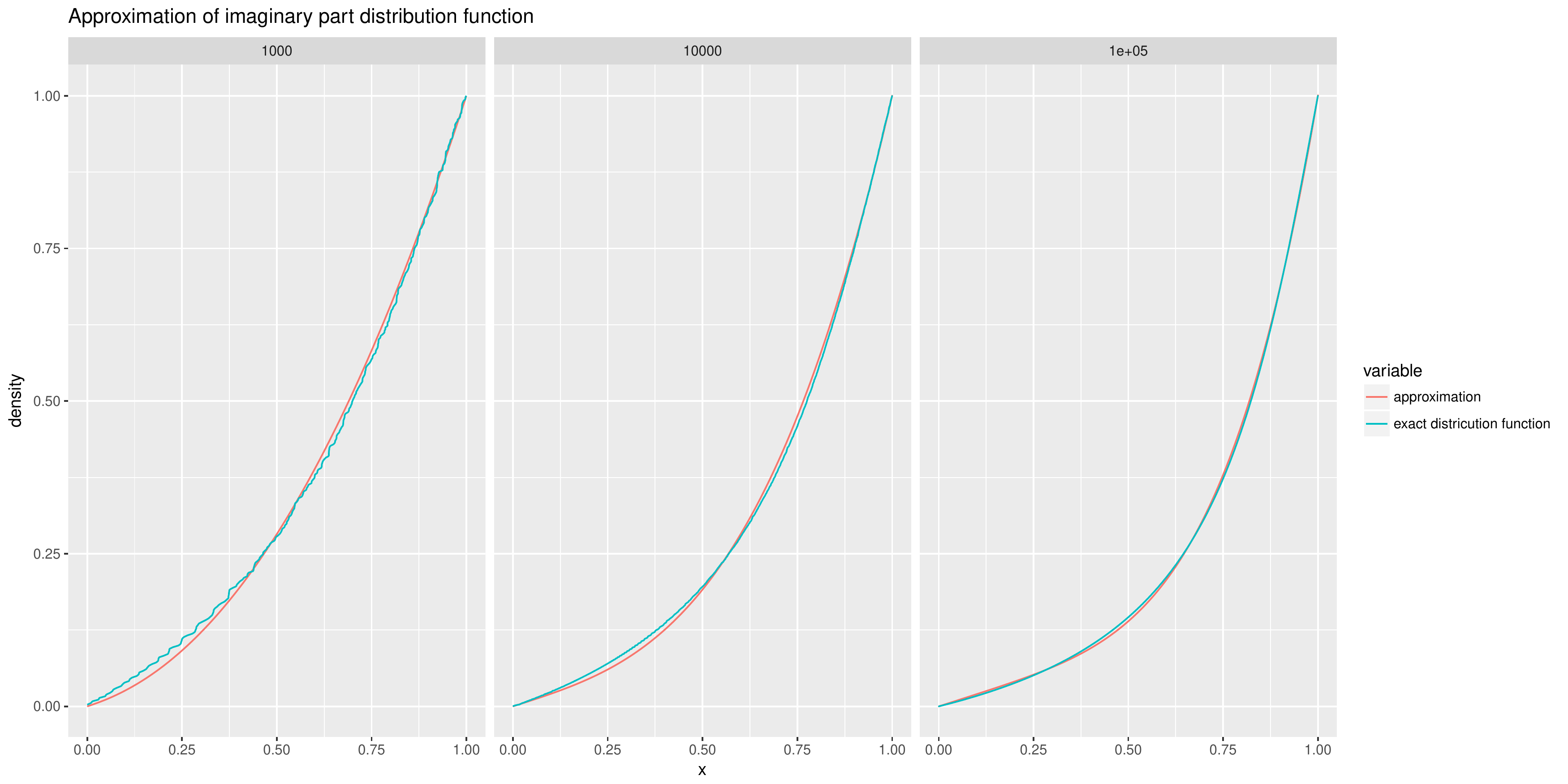}
	\caption{Approximation of the distribution function of the imaginary part for values inside the cluster at the vertex. $n$ was chosen to be 1000, 10000, and 100000, respectively.}
	\label{fig:im_distribution}
\end{figure}

Further, we analyse the distribution of the Fourier coefficients $z_j^{(n)}$ of the Riemann roots inside the cloud determined by the algorithm described above.
Therefore, we consider only the first half $2 \leq j \leq \frac{n}{2}$ as the coefficients are symmetric with respect to the real line ($z_j^{(n)} = \overline{z}_{n+1-j}^{(n)}$).
By normalizing the imaginary part to the interval $[0,1]$ we observe the following phenomena.
For greater $n \in \mathbb{N}$ is the more Fourier coefficients are clustered at the real line (which corresponds to a normalized value of 1).
\autoref{fig:im_distribution} shows how the probability function is getting more and more mass at larger values of the unit interval.
Another significant observation is that the probability function seem to have a slope of 0.2 at the point $x = 0$, regard less how the value $n$ is chosen.
In total the probability function can be approximated quite well with the ansatz function
$$ f(x) = 0.2 x + 0.8  x^p $$
where $p$ is chosen such that the integral fits the exact probability function. 

\begin{figure}[]
	\centering
	\includegraphics[width=\textwidth]{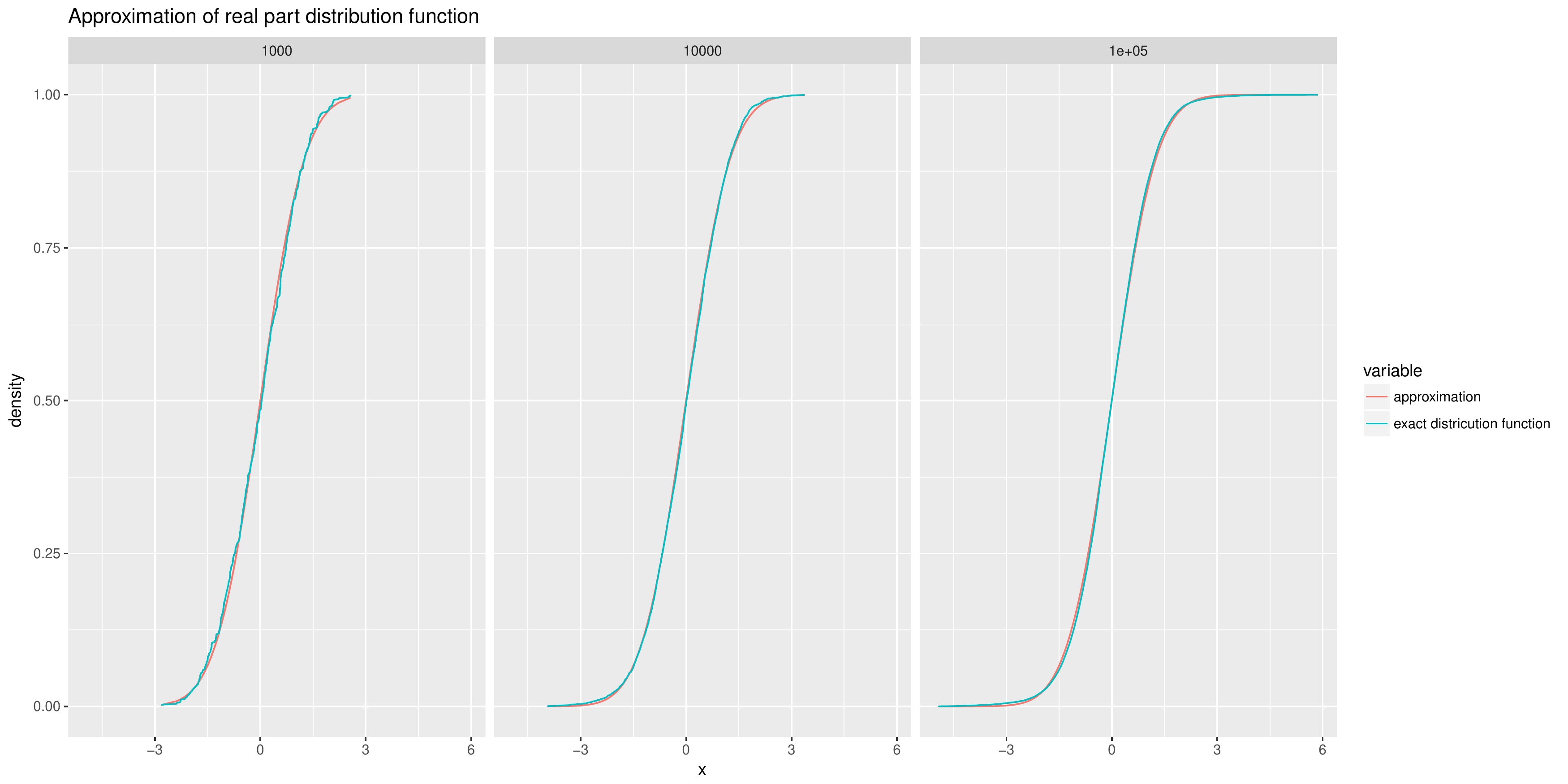}
	\caption{Approximation of the distribution function of the real part for values inside the cluster at the vertex. $n$ was chosen to be 1000, 10000, and 100000, respectively.}
	\label{fig:re_distribution}
\end{figure}

The approximation of the real part distribution function is more involved as the determination of the cloud index has quite a non-negative variance and the normalization to a unit interval of the real parts would be significantly dependent on $n$ and the cloud index algorithm.
Therefore, the normalization is not performed with the extreme values of $\Re(z)$, but rather with its variance.
In particular, we normalize the real part to have a mean 0 and variance 1.
We observe that the normalized real parts seem to be normally distributed, see \autoref{fig:re_distribution}.

\subsection{Nested Fourier analysis}

In the previous section the non-trivial roots of the Riemann zeta function were transformed by the discrete Fourier transform on the complex plane.
This open a new view at the distribution of these values.
In the following this approach is applied iteratively on the resulting numbers.

In particular, let $n \in \mathbb{N}$ be fixed and denote by $z_{(r)}$ and $z_{(i)}$ the vectors of real and imaginary parts of the Fourier coefficients of the Riemann roots $\gamma_j$ ($j = 1,\ldots,n$). These values are given by
\begin{align}
z_{(r)} &= \left(\sum_{j = 1}^n \frac{1}{\sqrt{n}} \gamma_j \cos\left(\frac{2\pi }{n}(j-1)(k-1)\right)\right)_{k = 1,\ldots,n} \\
z_{(i)} &= \left(\sum_{j = 1}^n \frac{1}{\sqrt{n}} \gamma_j \sin\left(\frac{2\pi }{n}(j-1)(k-1)\right)\right)_{k = 1,\ldots,n}
\end{align}

\begin{figure}[]
	\centering
	\includegraphics[width=\textwidth]{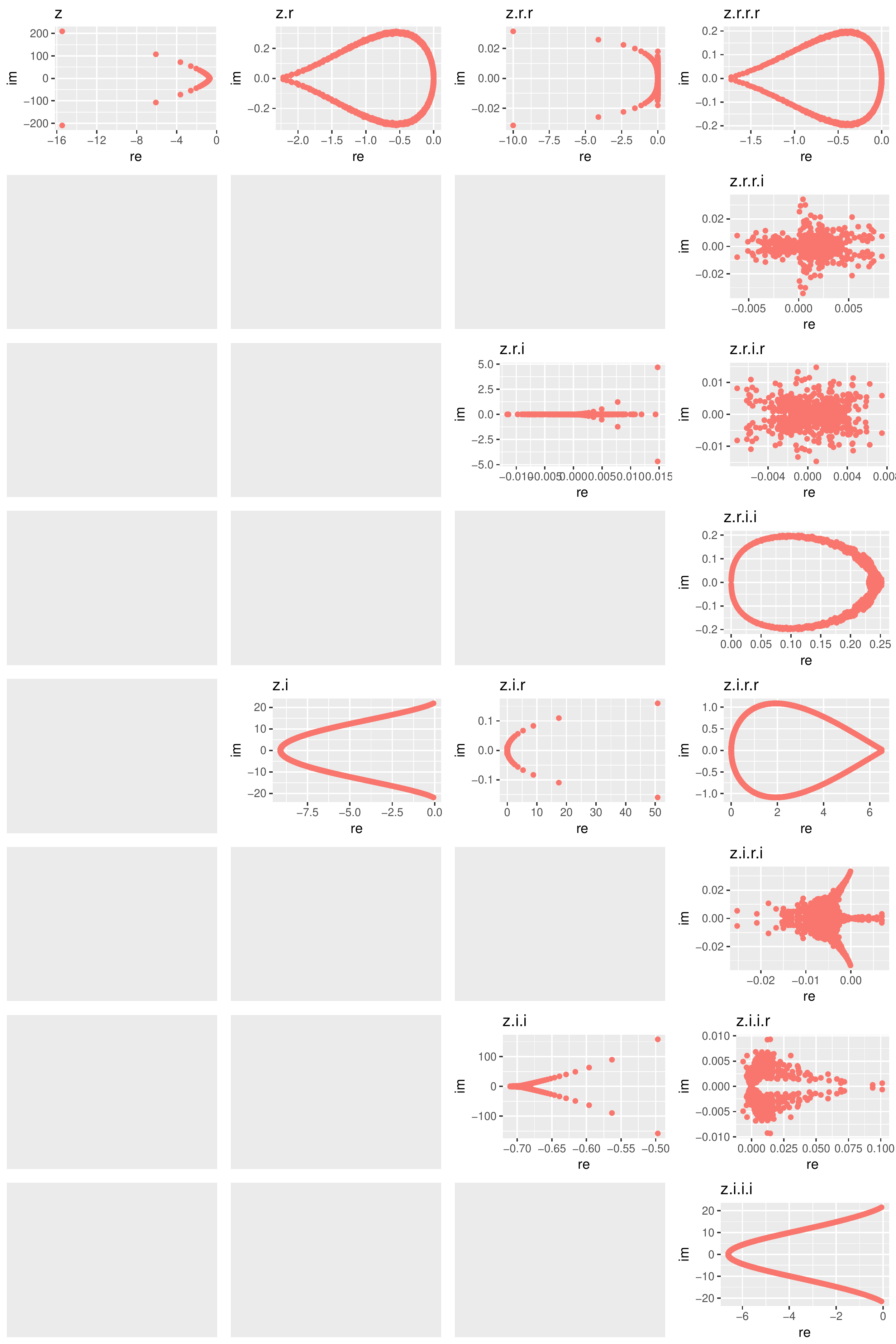}
	\caption{Nested Fourier analysis for $n = 1000$.}
	\label{fig:tree}
\end{figure}

Hereby, the first coefficient $(z_{(r)})_1$ equals $\sqrt{n}$ times the mean of $\gamma_k$ and $(z_{(i)})_1 = 0$.
In the following we neglect these two values, as they have only information about the average values and not the local distribution of the Riemann roots.
We rather apply the discrete Fourier transform on the vectors $(z_{(r)})_k$ and $(z_{(i)})_k$, where $k = 2,\ldots,n$.
Hence, we arrive at 
\begin{align}
z_{(x, r)} &= \left(\sum_{j = 1}^n \frac{1}{\sqrt{n-1}} (z_{(x)})_{j+1} \cos\left(\frac{2\pi }{n-1}(j-1)(k-1)\right)\right)_{k = 1,\ldots,n-1} \\
z_{(x, i)} &= \left(\sum_{j = 1}^n \frac{1}{\sqrt{n-1}} (z_{(x)})_{j+1} \sin\left(\frac{2\pi }{n-1}(j-1)(k-1)\right)\right)_{k = 1,\ldots,n-1}
\end{align}
where $x \in \{r, i\}$.
By iterating this procedure we end up with a binary tree of vectors
$z_{(x_1, \ldots, x_n)}$, where $x_l \in \{r, i\}$ for $l = 1,\ldots,n$.
this process stops a with each step the vectors get reduced in dimension by one and, hence, the tree has depth $n$ with $2^n$ nodes.
Together with the averages values that are neglected this tree inherits all information about the Riemann roots.
See \autoref{fig:tree} for the coefficients for $n = 1000$ and tree depth equal to 3.

\section*{Conclusion}
In this article we presented a new view on the Riemann hypothesis by means of the Fourier coefficients of the non-trivial roots of the Riemann zeta-function.
The transformation of these roots that form a real valued sequence on the complex numbers shows an interesting structure to study.
Together with the use of cyclic hermitian matrices defined here this can lead to a proof of the Hilbert-Polya conjecture.


\bibliographystyle{plain}
\bibliography{2018_Vartziotis_Merger_Contributions_to_the_study_of_the_non-trivial_roots_of_the_Riemann_zeta-function}


\end{document}